\pdfoutput=1
\RequirePackage{ifpdf}
\ifpdf %We are running pdfTeX in pdf mode
\documentclass[pdftex]{sigma}
\else
\documentclass{sigma}
\fi

\numberwithin{equation}{section}
\newtheorem{Theorem}{Theorem}[section]

\newtheorem{Lemma}[Theorem]{Lemma}
\newtheorem{Proposition}[Theorem]{Proposition}
{\theoremstyle{definition}
\newtheorem{Definition}[Theorem]{Definition}

\newtheorem{Example}[Theorem]{Example}
\newtheorem{Remark}[Theorem]{Remark} }

\DeclareMathOperator{\diag}{diag}
\DeclareMathOperator{\supp}{supp}
\DeclareMathOperator{\Ree}{Re}

\begin{document}

\newcommand{\arXivNumber}{1409.4213}

\allowdisplaybreaks

\renewcommand{\PaperNumber}{013}

\FirstPageHeading

\ShortArticleName{A Central Limit Theorem for Random Walks on the Dual of a~Compact Grassmannian}

\ArticleName{A Central Limit Theorem for Random Walks\\
on the Dual of a~Compact Grassmannian}

\Author{Margit R\"OSLER~$^\dag$ and Michael VOIT~$^\ddag$}

\AuthorNameForHeading{M.~R\"osler and M.~Voit}

\Address{$^\dag$~Institut f\"ur Mathematik, Universit\"at Paderborn,\\
\hphantom{$^\dag$}~Warburger Str.~100, D-33098 Paderborn, Germany}
\EmailD{\href{mailto:roesler@math.upb.de}{roesler@math.upb.de}}

\Address{$^\ddag$~Fakult\"at f\"ur Mathematik, Technische Universit\"at Dortmund,\\
\hphantom{$^\ddag$}~Vogelpothsweg 87, D-44221 Dortmund, Germany}
\EmailD{\href{mailto:michael.voit@math.uni-dortmund.de}{michael.voit@math.uni-dortmund.de}}

\ArticleDates{Received October 14, 2014, in f\/inal form February 03, 2015; Published online February 10, 2015}

\Abstract{We consider compact Grassmann manifolds $G/K$ over the real, complex or quaternionic numbers whose spherical
functions are Heckman--Opdam polynomials of type~$BC$.
From an explicit integral representation of these polynomials we deduce a~sharp Mehler--Heine formula, that is an
approximation of the Heckman--Opdam polynomials in terms of Bessel functions, with a~precise estimate on the error term.
This result is used to derive a~central limit theorem for random walks on the semi-lattice parametrizing the dual
of~$G/K$, which are constructed by successive decompositions of tensor powers of spherical representations of~$G$.
The limit is the distribution of a~Laguerre ensemble in random matrix theory.
Most results of this paper are established for a~larger continuous set of multiplicity parameters beyond the group
cases.}

\Keywords{Mehler--Heine formula; Heckman--Opdam polynomials; Grassmann manifolds; spherical functions; central limit
theorem; asymptotic representation theory}

\Classification{33C52; 43A90; 60F05; 60B15; 43A62; 33C80; 33C67}

\section{Introduction}

For Riemannian symmetric spaces $G/K$ of the compact or non-compact type, there is a~well-known contraction principle
which states that under suitable scaling, the spherical functions $\varphi_\lambda$ of $G/K$ tend to the spherical
functions $\psi_\lambda$ of the tangent space of $G/K$ in the base point, which is a~symmetric space of the f\/lat type:
\begin{gather*}
\lim_{n\to\infty} \varphi_{n\lambda}(\exp(x/n))=\psi_\lambda(x).
\end{gather*}
See~\cite{C} and, for a~more recent account,~\cite{BO}.
This curvature limit, also known as Mehler--Heine formula, extends to the more general setting of hypergeometric
functions associated with root systems, which converge under rescaling to generalized Bessel functions.
This is proven in~\cite{dJ} by a~limit transition in the Cherednik operators; see also~\cite{BO} for a~dif\/ferent
approach.
In the compact rank one case, the contraction principle is a~weak version of the classical Hilb formula for Jacobi
polynomials (see~\cite[Theorem 8.21.12]{Sz}), which provides in addition a~precise estimate on the rate of
convergence.
In this paper, we prove a~Mehler--Heine formula with a~precise estimate on the error term
for a~certain class of orthogonal polynomials associated with root systems, which in particular encompasses the
spherical functions of compact Grassmannians.
This result is a~``compact'' analogue of Theorem 5.4 in~\cite{RV2}, which gives a~scaling limit with error bounds for
hypergeometric functions in the dual, non-compact setting.
In the second part of the paper, we shall use the Mehler--Heine formula~\ref{vgl-Bessel} in order to establish a~central
limit theorem for random walks on the semi-lattice of dominant weights parametrizing the unitary dual of a~compact
Grassmannian.

To become more precise, we consider the compact Grassmannians $\mathcal G_{p,q}(\mathbb F)=G/K$ over one of the (skew-)
f\/ields $\mathbb F= \mathbb R, \mathbb C, \mathbb H$, with $G=SU(p+q,\mathbb F)$ and $K=S(U(q, \mathbb F)\times
U(p,\mathbb F))$, where $p\ge q\ge 1$.
Via polar decomposition of~$G$, the double coset space $G//K=\{KgK:\> g\in G\}$ may be topologically identif\/ied with the
fundamental alcove
\begin{gather*}
A_0:=\Big\{x=(x_1,\ldots,x_q)\in\mathbb R^q:\> \frac{\pi}{2}\ge x_1\ge x_2\ge\dots\ge x_q\ge0\Big\},
\end{gather*}
with $x\in A_0$ being identif\/ied with the matrix
\begin{gather*}
a_x=
\begin{pmatrix}
\cos\underline x & -\sin \underline x & 0
\\
\sin \underline x & \cos \underline x & 0
\\
0 & 0 & I_{p-q}
\end{pmatrix}
.
\end{gather*}
Here we use the diagonal matrix notation $\underline x=\diag(x_1,\ldots,x_q)$, and the functions $\sin$, $\cos$ are
understood component-wise.
For details, see~\cite[Theorem~4.1]{RR}.
The spherical functions of~$\mathcal G_{p,q}(\mathbb F)$ can be viewed as Heckman--Opdam polynomials of type~$BC_q$,
which are also known as multivariable Jacobi polynomials.
They may be described as follows: denote by~$F_{BC}(\lambda,k;\cdot)$ the Heckman--Opdam hypergeometric function
associated with the root system
\begin{gather*}
R=2  BC_q=\{\pm 2e_i, \pm 4 e_i: \> 1\leq i\le q\}\cup\{\pm 2e_i \pm 2e_j: 1\leq i < j \leq q\} \subset \mathbb R^q,
\end{gather*}
with spectral variable $\lambda\in\mathbb C^q$ and multiplicity parameter $k=(k_1, k_2, k_3)\in \mathbb R^3$
corresponding to the roots $\pm 2 e_i$, $\pm 4 e_i$ and $2(\pm e_i\pm e_j)$.
Fix the positive subsystem
\begin{gather*}
R_+=\{2e_i, 4e_i, 1\leq i \leq q\}\cup\{2e_i\pm 2e_j, 1\leq i<j \leq q\}
\end{gather*}
and the associated semi-lattice of dominant weights,
\begin{gather*}
P_+=\big\{\lambda\in (2\mathbb Z)^q: \> \lambda_1\ge\lambda_2\ge\dots\ge\lambda_q\ge0\big\}.
\end{gather*}
Then the set of spherical functions of $\mathcal G_{p,q}(\mathbb F)$ is parametrized by $P_+$ and consists of the
functions
\begin{gather}
\label{jacobi-pol-def}
\varphi_\lambda^p(a_x)=F_{BC}(\lambda+\rho_p, k(p),ix)=: R_\lambda^p(x),
\qquad
\lambda\in P_+
\end{gather}
with multiplicity parameter
\begin{gather}
\label{multipl}
k(p)=(d(p-q)/2, (d-1)/2, d/2),
\end{gather}
where $d=\dim_\mathbb R \mathbb F \in \{1,2,4\}$ and
\begin{gather*}
\rho_p=\frac{1}{2}\sum\limits_{\alpha\in R_+} k_\alpha\alpha=\sum\limits_{i=1}^q \left(\frac{d}{2}(p+q+2-2i) -1\right)e_i.
\end{gather*}
The functions $R_\lambda^p$ are the Heckman--Opdam polynomials associated with the root system~$R$ (called Jacobi
polynomials in the following) and with multiplicity $k(p)$, normalized according to $R_\lambda^p(0)=1$.
We refer to~\cite{H,HS,O1} for Heckman--Opdam theory in general, and to~\cite{RR} and the references
cited there for the connection with spherical functions in the compact $BC$ case.
Notice that our notion of~$F_{BC}$ coincides with that of Heckman, Opdam and~\cite{R2,RV2}, while it dif\/fers
from the geometric notion in~\cite{RR}.
Theorem~4.3 of~\cite{RR} corresponds to~\eqref{jacobi-pol-def}.

In Theorem 4.2 of~\cite{RR}, the product formula for spherical functions of $(G,K)$ was written as a~formula on~$A_0$
and analytically extended to a~product formula for the Jacobi polyno\-mials~$R_\lambda^p$ with multiplicity~$k(p)$
corresponding to arbitrary real parameters $p>2q-1$.
This led to three continuous series of positive product formulas for Jacobi polynomials corresponding to $\mathbb
F=\mathbb R, \mathbb C, \mathbb H$ and to associated commutative hypergroup structures on $A_0$; see~\cite{J}
and~\cite{BH} for the notion of hypergroups.
Using a~Harish-Chandra-type integral representation for the $R_\lambda^p$, we shall derive a~Mehler--Heine formula with
a~precise asymptotic estimate for the Jacobi polyno\-mials~$R_{\lambda}^p$ in terms of Bessel functions associated with
root system $B_q$ on the Weyl chamber
\begin{gather*}
C=\{x= (x_1, \ldots, x_q) \in \mathbb R^q: x_1 \geq \cdots \geq x_q \geq 0\}.
\end{gather*}
This Mehler--Heine formula will be the key ingredient for the main result of the present paper, a~central limit theorem
for random walks on the semi-lattice $P_+$, which parametrizes the spherical unitary dual of $G/K$.
To explain this CLT, let us f\/irst recall that via the GNS representation, the spherical functions $\varphi_\lambda,
\lambda \in P_+$ of $(G,K)$, which are necessarily positive def\/inite, are in a~one-to-one correspondence with the
(equivalence classes of) spherical representations $(\pi_\lambda,H_\lambda)$ of~$G$, that is those irreducible unitary
representations of~$G$ whose restriction to~$K$ contains the trivial representation with multiplicity one, see~\cite{F}
or~\cite[Chapter~IV]{Hel}.
The decomposition of tensor products of spherical representations into their irreducible components leads to
a~probability preserving convolution~$*_{d,p}$ and f\/inally a~Hermitian hypergroup structure on the discrete set~$P_+$;
see~\cite{Du} and~\cite{K3}.
Following, e.g.,~\cite{BH,V1,Z}, we introduce random walks $(S_n^{d,p})_{n\ge0}$ on $P_+$ associated
with $*_{d,p}$ and derive some limit theorems for $n\to\infty$.
The main result of this paper will be the Central Limit Theorem~\ref{clt}.
This CLT implies the following result for $\mathcal G_{p,q}(\mathbb F)=G/K$:

\begin{Theorem}
\label{clt-spezial}
Let $(\pi_\lambda,H_\lambda)$ be a~non-trivial spherical representation of~$G$ associated with $\lambda\in
P_+\setminus\{0\}$.
Let $u_\lambda\in H_\lambda$ be~$K$-invariant with $\|u_\lambda\|=1$.
For each $n\in \mathbb N$, decompose the~$n$-fold tensor power $(\pi_\lambda^{\otimes, n},H_\lambda^{\otimes, n})$ into
its finitely many irreducible unitary components
\begin{gather*}
\big(\pi_\lambda^{\otimes, n},H_\lambda^{\otimes, n}\big)=\left(\bigoplus_{\tau_n}\pi_{\tau_n}, \bigoplus_{\tau_n}H_{\tau_n}\right),
\end{gather*}
where the components are counted with multiplicities.
Consider the orthogonal projections $p_{\tau_n}:H_\lambda^{\otimes, n}\to H_{\tau_n}$ and a~$P_+$-valued random variable
$X_{n,\lambda}$ with the finitely supported distribution
\begin{gather*}
\sum\limits_{\tau_n}\big\|p_{\tau_n}(u_\lambda^{\otimes, n})\big\|^2 \delta_{\tau_n}\in M^1(P_+)
\end{gather*}
with the point measures $\delta_{\tau_n}$ at $\tau_n$.
Then, for $n\to\infty$,
\begin{gather*}
\frac{X_{n,\lambda}}{m(\lambda) \sqrt n}
\end{gather*}
tends in distribution to
\begin{gather*}
d\rho_{d,p}(x)= c_{d,p}^{-1} \prod\limits_{j=1}^q x_j^{d(p-q+1)-1}  \prod\limits_{1\le i<j\le q}\big(x_i^2-x_j^2\big)^d
  e^{-(x_1^2+\cdots+x_q^2)/2}\> dx
\in M^1(C)
\end{gather*}
with a~suitable normalization $ c_{d,p}$.
Notice that the probability measure $d\rho_{d,p}$ is the distribution of a~Laguerre ensemble on~$C$.
The modified variance parameter $m(\lambda)>0$ is a~second order polynomial in~$\lambda$ and given explicitly in
Lemma~{\rm \ref{basic-prop-moment}} below.
\end{Theorem}

For $q=1$, the Central Limit Theorem~\ref{clt} has a~long history as a~CLT for random walks on~$\mathbb Z_+$ whose
transition probabilities are related to product linearizations of Jacobi polynomials.
This includes random walks on the duals of ${\rm SU}(2)$ and $({\rm SO}(n),{\rm SO}(n-1))$ in~\cite{ER} and~\cite{Ga}.
See also~\cite{V1} for further one-dimensional cases.
For $q\ge2$ our results are very closely related to the work~\cite{CR} of Clerc and Roynette on duals of compact
symmetric spaces.
For a~survey on limits for spherical functions and CLTs in the non-compact case for $q=1$ we refer to~\cite{V2}.

\section{A Mehler--Heine formula}\label{Section2}

In this section we derive a~Mehler--Heine formula for the Jacobi polynomials $R_\lambda^p (\lambda \in P_+)$, describing
the approximation of these polynomials in terms of Bessel functions with a~precise error bound.
Our result will be based on Laplace-type integrals for the Jacobi polynomials and the associated Bessel functions, where
we treat the group cases with integers $p\ge q$ as well as the case $p\in\mathbb R$ with $p\ge 2q-1$ beyond the group
case.
The integral representation for $R_{\lambda}^p$ below is a~special case of a~more general Harish-Chandra integral
representation for hypergeometric functions $F_{BC}$ in~\cite{RV2}.
To start with, let us introduce some notation:

Let $H_q(\mathbb F)=\{x\in M_q(\mathbb F): x^*:= \overline x^t=x\}$ denote the space of Hermitian matrices over
$\mathbb F$, and denote by $\Delta(x)$ the determinant of $x\in H_q(\mathbb F)$, which may be def\/ined as the product of
(right) eigenvalues of~$x$.
We mention that for $\mathbb F=\mathbb H$, this is just the Moore determinant, which coincides with the Dieudonn\'e
determinant if~$x$ is positive semi-def\/inite, see, e.g.,~\cite{A}.
On $H_q(\mathbb F)$, we consider the power functions
\begin{gather*}
\Delta_\lambda(a):= \Delta_1(a)^{\lambda_1-\lambda_2}   \cdots   \Delta_{q-1}(a)^{\lambda_{q-1}-\lambda_q} \Delta_q(a)^{\lambda_{q}},
\qquad
\lambda \in \mathbb C^q,
\end{gather*}
with the principal minors $\Delta_r(a)=\det ((a_{ij})_{1\leq i,j\leq r}) $ of the matrix $a=(a_{ij})_{1\leq i,j\leq
q}\in H_q(\mathbb F)$, see~\cite{FK}.
We introduce the matrix ball $ B_q:= \{w\in M_q(\mathbb F): w^*w < I\}$, where $A<B$ means for matrices $A,B\in
M_q(\mathbb F)$ that $B-A$ is (strictly) positive def\/inite.
On $B_q$, we def\/ine the probability measures
\begin{gather*}
dm_p(w)=\frac{1}{\kappa_{pd/2}} \Delta(I-w^*w)^{pd/2-\gamma}dw
\in M^1(B_q),
\end{gather*}
with $p\in \mathbb R$, $p>2q-1$.
Here $dw$ is the Lebesgue measure on the ball $B_q$,
\begin{gather*}
\gamma:= d\left(q-\frac{1}{2}\right)+1
\end{gather*}
and
\begin{gather*}
\kappa_{pd/2}= \int_{B_q} \Delta(I-w^*w)^{pd/2-\gamma}\> dw.
\end{gather*}
According to Theorem 2.4 of~\cite{RV2}, the Heckman--Opdam hypergeometric function $F_{BC}(\lambda, k(p),x)$ with
$\lambda \in \mathbb C^q$, $x\in \mathbb R^q$ and $k(p)$ as in~\eqref{multipl} has the following integral representation
for $p\in \mathbb R$ with $p>2q-1$:
\begin{gather*}
F_{BC}(\lambda, k(p), x)=\int_{B_q\times U_0(q, \mathbb F)} \Delta_{(\lambda-\rho_p)/2}(g_x(u,w))\> dm_p(w)du,
\end{gather*}
where $U_0(q,\mathbb F)$ denotes the identity component of $U(q, \mathbb F)$ and
\begin{gather*}
g_x(u,w)=u^{-1}(\cosh \underline x+w^* \sinh \underline x)(\cosh \underline x+\sinh\underline x w)u.
\end{gather*}
It is easily checked that $U_0(q, \mathbb F)$ may be replaced by $U(q,\mathbb F)$ in the domain of integration.
Notice further that $x\mapsto g_x(u,w)$ extends to a~holomorphic function on $\mathbb C^q$.
As the principal minors~$\Delta_r(a)$ are polynomial in the entries of $a\in H_q(\mathbb F)$, it follows
that $x\mapsto\Delta_{\lambda/2}(g_x(u,w))$ extends to a~holomorphic function on~$\mathbb C^q$ for each~$\lambda \in P_+$.
In view of relation~\eqref{jacobi-pol-def}, this leads to the following integral representation for the Jacobi
polynomials~$R_\lambda^p$:

\begin{Proposition}
\label{int-rep-prop}
Let $p\in \mathbb R$ with $p>2q-1$ and $k(p)= (d(p-q)/2, (d-1)/2, d/2)$ with $d\in\{1,2,4\}$.
Then the Jacobi polynomials $ R_\lambda^p$, $\lambda \in P_+$, have the integral representation
\begin{gather}
\label{int-rep}
R_\lambda^p(x)= \int_{B_q\times U(q,\mathbb F)} \Delta_{\lambda/2}(g_{ix}(u,w))\> dm_p(w)du
\qquad
\text{for}
\quad
x\in A_0
\end{gather}
with
\begin{gather*}
g_{ix}(u,w)=u^{-1}(\cos\underline x+w^* i\sin \underline x)(\cos\underline x+i\sin \underline x w)u.
\end{gather*}
\end{Proposition}

We next turn to the Bessel functions which will show up in the Mehler--Heine formula.
They are given in terms of Bessel functions of Dunkl type which generalize the spherical functions of Cartan motion
groups; see~\cite{dJ} and~\cite{O} for a~general background.
We denote by $J_k^B$ the Bessel function which is associated with the rational Dunkl operators for the root system $B_q=
\{\pm e_i, \pm e_i \pm e_j: 1\leq i< j \leq q\}$ and multiplicity $k=(k_1, k_2)$ corresponding to the roots $\pm e_i$
and $\pm e_i \pm e_j$.
We shall be concerned with multiplicities which are connected as follows to the $BC_q$ multiplicities $k(p)$
from~\eqref{multipl}:
\begin{gather*}
k= (k_1, k_2)
\qquad
\text{with}
\quad
k_1=k(p)_1+k(p)_2=d(p-q+1)/2-1/2,
\quad
k_2=k(p)_3=d/2.
\end{gather*}
For such~$k$ on $B_q$, we use the notion
\begin{gather*}
\widetilde\varphi_\lambda^{p}(x):=J_k^B(x, i\lambda),
\qquad
x\in C, \lambda \in \mathbb C^q.
\end{gather*}
It is well-known that for integers $p\geq q$, the $\widetilde\varphi_\lambda^{p}$ are the spherical functions of the
Euclidean symmetric spaces $G_0/K$, where $K=S(U(p,\mathbb F)\times U(q,\mathbb F))$ and $G_0= K\ltimes M_{p,q}(\mathbb F)$
is the Cartan motion group associated with the Grassmannian $\mathcal G_{p,q}(\mathbb F)$. Hereby the double coset
space $G_0//K$ is identif\/ied with the Weyl chamber~$C$ such that $x\in C$ corresponds to the double coset of $(I_p\times
I_q, (I_{p-q},\underline x) ) \in G_0$, and in this way,~$K$-biinvariant functions on $G_0$ may be considered as
functions on~$C$.
Two functions $\widetilde\varphi_\lambda^{p}$ and $\widetilde\varphi_\mu^{p}$ coincide if and only if~$\lambda$
and~$\mu$ are in the same Weyl group orbit.
Finally, the bounded spherical functions are exactly those $\widetilde\varphi_\lambda^{p}$ with $\lambda\in C$.
The Bessel functions $\widetilde\varphi_\lambda^{p}$ with $d=\dim_\mathbb F\mathbb R$ and not necessarily
integral parameter~$p$ are closely related to Bessel functions on the symmetric cone of positive def\/inite $q\times
q$-matrices over $\mathbb F$, see Section~4 of~\cite{R1}.
It has been shown there that for $p>2q-1$, they have a~positive product formula which generalizes the product formula in
the Cartan motion group cases and leads to a~commutative hypergroup structure on the Weyl chamber~$C$.

\begin{Lemma}
\label{Bessel-inte}
For $p\in \mathbb R$ with $p>2q-1$, the Bessel functions $\widetilde\varphi_\lambda^{p}$ with $\lambda \in \mathbb R^q$
have the following integral representation:
\begin{gather}
\label{int-rep-Bessel}
\widetilde\varphi_\lambda^{p}(x)=\int_{B_q}\int_{U(q,\mathbb F)} e^{i\Ree \operatorname{tr}(w \underline x u \underline\lambda)}dm_p(w)du.
\end{gather}
\end{Lemma}

\begin{proof}
This follows readily from equations~(3.12) and~(4.4) in~\cite{R1}; see also Proposition 5.3 of~\cite{RV2}.
\end{proof}

\begin{Remark}
\label{group}
There are f\/initely many geometric cases which are not covered by the range $p\in{}]2q-1,\infty[$, namely the indices
$p\in\{q, q+1, \ldots, 2q-1\}$. In these cases, the Jacobi polyno\-mials~$R_\lambda^p $ and the Bessel functions
$\widetilde\varphi_\lambda^{p}$ both admit interpretations as spherical functions and have an integral representation
similar to that above, by the following reasoning: According to Lemma~2.1 of~\cite{R2}, the measure $m_p\in M^1(B_q)$
with $p\in \mathbb N$, $p\geq 2q$ is just the pushforward measure of the normalized Haar measure on $U(p,\mathbb F)$ under
the mapping
\begin{gather*}
v\mapsto \sigma_0^* v\sigma_0,
\qquad
\text{with}
\quad
\sigma_0=
\begin{pmatrix}
I_q
\\
0_{(p-q)\times q}
\end{pmatrix}
\in M_{p,q}(\mathbb F).
\end{gather*}
For $p\in \{q, q+1, \ldots, 2q-1\}$, we now def\/ine the measure $m_p\in M^1(B_q)$ in the same way as a~pushforward
measure of the Haar measure on $U(p, \mathbb F)$. (But in contrast to the case $p\geq 2q$, it will not have a~Lebesgue
density in these cases).
From the integral representa\-tions~(3.3) and~(4.4) of~\cite{R1} for the Bessel functions, as well as Theorem~2.1
of~\cite{RV2} and relation~\eqref{jacobi-pol-def} between Jacobi polynomials and hypergeometric functions, one obtains
that the integral representations of Proposition~\ref{int-rep-prop} and Lemma~\ref{Bessel-inte} extend
to the case $p\in\{q, q+1, \ldots, 2p-1\}$.
\end{Remark}

We shall now compare the integral representations of Proposition~\ref{int-rep-prop} and Lemma~\ref{Bessel-inte}, which
will lead to the following quantitative Mehler--Heine (or Hilb-type) formula.

\begin{Theorem}
\label{vgl-Bessel}
There exist constants $C_1,C_2>0$ such that for all $p\in\{q,q+1,\ldots,2q-1\}\cup$ $] 2q-1,\infty[$, all $\lambda\in
P_+$, and $ x\in A_0$,
\begin{gather*}
\left|R_{\lambda}^p(x)-\widetilde\varphi_\lambda^{p} (x)\right|\leq C_1  x_1^2\lambda_1  e^{C_2 x_1^2\lambda_1}.
\end{gather*}
\end{Theorem}
Thus in particular,
\begin{gather*}
\left|R_{n\lambda}^p \left(\frac{x}{n}\right)-\widetilde\varphi_\lambda^{p}(x)\right|\leq \frac{C_1}{n}
x_1^2\lambda_1  e^{C_2 x_1^2\lambda_1/n} \to 0
\qquad
\text{for}
\quad
n\to \infty.
\end{gather*}
Notice that the estimate of Theorem~\ref{vgl-Bessel} is uniform in~$p$, a~fact which was to our knowledge so far not
even noticed in the rank-one case.
We conjecture that the statement of this theorem remains correct for $p\in [q,\infty[$.
\begin{proof}
We only consider the case $p>2q-1$ where the proof is based on Proposition~\ref{int-rep-prop} and
Lemma~\ref{Bessel-inte}.
By the previous remark, the cases $p={q,q+1,\ldots,2q-1}$ can be treated in the same way.
Notice that it suf\/f\/ices to check uniformity in~$p$ for $p>2q-1$.

We substitute $w\mapsto u^*w^*$ in the integral~\eqref{int-rep-Bessel} and obtain
\begin{gather*}
\widetilde\varphi_\lambda^{p}(x)= \int_{B_q\times U(q, \mathbb F)} e^{i \Ree \operatorname{tr}(u^* w^*\underline x u\underline \lambda )}dm_p(w)du.
\end{gather*}
Denoting the trace of the upper left $(r\times r)$-block of a~$(q\times q)$-matrix by $\operatorname{tr}_r$, we have
\begin{gather*}
\Ree \operatorname{tr}(u^* w^*\underline x u\underline\lambda)
=\frac{1}{2} \sum\limits_{r=1}^q (u^*((\underline x w)^*+\underline x w)u)_{rr}  \lambda_r
\\
\phantom{\Ree \operatorname{tr}(u^* w^*\underline x u\underline\lambda)}{}
=\sum\limits_{r=1}^q \bigl[ \operatorname{tr}_r(u^*((\underline xw)^*+\underline xw)u) - \operatorname{tr}_{r-1}(u^*((\underline xw)^*+\underline xw)u)\bigr]   \lambda_r/2
\\
\phantom{\Ree \operatorname{tr}(u^* w^*\underline x u\underline\lambda)}{}
=\sum\limits_{r=1}^q \operatorname{tr}_r(u^*((\underline xw)^*+ \underline xw)u)  (\lambda_r-\lambda_{r+1})/2
\end{gather*}
with $\lambda_{q+1}:=0$.
Hence
\begin{gather*}
\widetilde\varphi_\lambda^{p}(x)
=\int_{B_q\times U(q,\mathbb F)} \prod\limits_{r=1}^q e^{i  \operatorname{tr}_r(u^*((xw)^*+xw)u) (\lambda_r-\lambda_{r+1})/2} dm_p(w)du.
\end{gather*}
Furthermore, by Proposition~\ref{int-rep-prop},
\begin{gather*}
R_{\lambda}^p(x)=\int_{B_q\times U(q,\mathbb F)} \prod\limits_{r=1}^q\Delta_r(g_{ix}(u,w))^{(\lambda_r-\lambda_{r+1})/2} dm_p(w)du.
\end{gather*}
Telescope summation yields the well-known estimate
\begin{gather*}
\left|\prod\limits_{r=1}^q a_r
-\prod\limits_{r=1}^q b_r\right|\le\big(\max(|a_1|,\ldots,|a_q|,|b_1|,\ldots, |b_q|)\big)^{q-1}  \sum\limits_{r=1}^q|a_r-b_r|
\end{gather*}
for $a_1,\ldots,a_q,b_1,\ldots,b_q\in \mathbb C$.
We thus obtain
\begin{gather}
\big|R_{\lambda}^p(x)-\widetilde\varphi_\lambda^{p}(t)\big| \leq \sum\limits_{r=1}^q \int_{B_q\times U(q,\mathbb F)}
M(x,u,w,\lambda)
\label{main-estimate}
\\
\hphantom{\big|R_{\lambda}^p(x)-\widetilde\varphi_\lambda^{p}(t)\big| \leq}{}
\times \big| \Delta_r(g_{x}(u,w))^{(\lambda_r-\lambda_{r+1})/2}
 - e^{i \operatorname{tr}_r(u^*((\underline tw)^*+ \underline tw)u)  (\lambda_r-\lambda_{r+1})/2}\big| dm_p(w)du
\nonumber
\end{gather}
with
\begin{gather*}
M(x,u,w,\lambda):=\max \Big(1, \max_{r=1,\ldots,q}\big|\Delta_r(g_{x}(u,w))^{(\lambda_r-\lambda_{r+1})/2}\big|^{q-1}\Big).
\end{gather*}
We now investigate $\Delta_r(g_{ix}(u,w))^{(\lambda_r-\lambda_{r+1})/2}$ more closely.
As $x$, $u$, $w$ run through compacta, we obtain that uniformly in $x$, $u$, $w$,
\begin{gather*}
g_{ix}(u,w)=u^{-1}(\cos\underline x+w^*i \sin \underline x)(\cos\underline x+i\sin \underline x w)u
\\
\phantom{g_{ix}(u,w)}
=u^{-1}\big(I_q+ w^*i\underline x+O\big(x^2\big)\big) \big(I_q+i \underline x w+O\big(x^2\big)\big)u
\\
\phantom{g_{ix}(u,w)}
=I_q+ u^{-1}(i\underline x w+w^*i\underline x)u +O\big(x^2\big),
\end{gather*}
and thus
\begin{gather}
\label{entwickel}
\Delta_r(g_{x}(u,w))= 1+ \operatorname{tr}_r\big(u^{-1}(i\underline xw+w^*i\underline x)u\big)+O\big(x^2\big).
\end{gather}
Using the power series for $\ln(1+z)$, we further have
\begin{gather*}
\Delta_r(g_{ix}(u,w) )^{(\lambda_r-\lambda_{r+1})/2}=\exp\big[\tfrac 12 (\lambda_r-\lambda_{r+1})  \ln\big(1+
\operatorname{tr}_r\big(u^{-1}(i\underline xw+w^*i\underline x)u\big) +O\big(x^2\big)\big)\big]
\\
\hphantom{\Delta_r(g_{ix}(u,w) )^{(\lambda_r-\lambda_{r+1})/2}}{}
 =\exp\big[\tfrac 12(\lambda_r-\lambda_{r+1})  \operatorname{tr}_r\big(u^{-1}(i\underline xw+w^*i\underline x)u\big)\!+O\big(x^2\big) (\lambda_r-\lambda_{r+1})\big].\!
\end{gather*}
Notice that $y:= u^{-1}(i\underline xw+w^*i\underline x)u $ is skew-Hermitian, that is $y^*=-y$.
Therefore $\overline {\operatorname{tr}_r(y)}=- \operatorname{tr}_r(y)$, which implies that $\Ree(\operatorname{tr}_r(y))=0$. It follows that
\begin{gather*}
\big|\Delta_r(g_{ix}(u,w))^{(\lambda_r-\lambda_{r+1})/2}\big|=\exp\big[\tfrac 12 (\lambda_r-\lambda_{r+1})
 \Ree\big(\operatorname{tr}_r(y)+ O\big(x^2\big)\big)\big]=e^{(\lambda_r -\lambda_{r+1}) O(x^2)}.
\end{gather*}
Note that these considerations apply for all f\/ields $\mathbb F=\mathbb R, \mathbb C, \mathbb H$. It follows that there
exists a~constant $C_3>0$ (independent of $x$, $u$, $w$, $\lambda$) such that
\begin{gather*}
M(x,u,w,\lambda)\le e^{C_3x_1^2\lambda_1}
\qquad
\text{for all}
\quad
x\in A_0, \  \lambda\in P_+,  \ u\in U(q, \mathbb F), \ w \in B_q.
\end{gather*}
From this inequality we obtain by the mean value theorem that for all $x\in A_0$ and $\lambda\in P_+$,
\begin{gather*}
\big|  \Delta_r(g_{x}(u,w))^{(\lambda_r-\lambda_{r+1})/2}- e^{i \operatorname{tr}_r(u^*((\underline x w)^*+ \underline x w)u)
 (\lambda_r-\lambda_{r+1})/2}\big|
\\
\qquad{}
 \leq e^{C_3x_1^2\lambda_1}-1 \leq C_3 x_1^2\lambda_1 e^{C_3x_1^2\lambda_1}.
\end{gather*}
These estimates together with~\eqref{main-estimate} imply the assertion.
\end{proof}

\begin{Example}[the rank one case] \sloppy
\label{1-dim-bsp-1}
For $q=1$ the Jacobi polyno\-mials~$R_\lambda^p$ associated with root system $BC_1=\{\pm e_1,\pm 2e_1\}$ are classical
one-variable Jacobi polynomials in trigonometric parametrization.
For integers~$p$, the associated Grassmannians are the projective spaces~$P_p(\mathbb F)$.
For the details, recall that the classical Jacobi polynomials $R_n^{(\alpha,\beta)}$ with the normalization
\mbox{$R_n^{(\alpha,\beta)}(1)=1$} are given~by
\begin{gather}
\label{class-Jacobi1}
R_n^{(\alpha,\beta)}(x)={}_2F_1(\alpha+\beta+n+1, -n;\alpha+1; (1-x)/2),
\end{gather}
where $n\in\mathbb Z_+$, $\alpha,\beta>-1$.
It is easily derived from the example on p.
17 of~\cite{O1} that
\begin{gather}
\label{class-Jacobi2}
R_\lambda^p(x)= R_{\lambda/2}^{(\alpha,\beta)} (\cos 2x)
\end{gather}
for $\lambda\in 2\mathbb Z_+$, with
\begin{gather*}
\alpha=(dp-2)/2,
\qquad
\beta=(d-2)/2;
\end{gather*}
see also~\cite[Section~5]{RR}.
In the rank one case, the $U(1, \mathbb F)$ integral in representation~\eqref{int-rep} cancels by invariance of~$\Delta$
under unitary conjugation.
Thus~\eqref{int-rep} reduces to
\begin{gather*}
R_\lambda^p(x)=\frac{1}{\kappa_{pd/2}}\int_{B_1} \left((\cos x+\overline w i \sin x)(\cos x+i\sin x
w)\right)^{\lambda/2}  \big(1-|w|^2\big)^{d(p-1)/2 -1} \> dw
\end{gather*}
for $\lambda \in \mathbb Z_+$, $p>1$. In particular, if $\mathbb F=\mathbb R$, then $d=1$ and $ B_1=[-1,1]$. Thus
\begin{gather*}
R_n^{(p/2-1, -1/2)}(\cos 2x)=\frac{1}{\kappa_{p/2}}\int_{-1}^1 (\cos x+it\sin x)^{2n} \big(1-t^2\big)^{(p-3)/2} \> dt.
\end{gather*}
If $\mathbb F=\mathbb C$, then $d=2$ and $B_1=\{z\in\mathbb C:\> |z|\le 1\}$.
Using polar coordinates $z=te^{i\theta}$, one obtains
\begin{gather*}
R_n^{(p-1, 0)}(\cos 2x)=\frac{1}{\kappa_{p}}\int_{-1}^1\int_{0}^\pi\! \big((\cos x+i te^{i\theta} \sin x)(\cos x+i
te^{-i\theta} \sin x)\big)^{2n} \big(1-t^2\big)^{p-2}t \> dt d\theta.
\end{gather*}
The quaternionic case can be treated in a~similar way.
These formulas are just special cases of a~well-known Laplace-type integral representation for Jacobi polynomials with
general indices $\alpha\ge\beta\ge -1/2$; see, e.g.,~\cite[Section~18.10]{OLBC}.
\end{Example}

{\sloppy Let us f\/inally mention that the Mehler--Heine formula~\ref{vgl-Bessel} corresponds to~\cite[Theo\-rem~8.21.12]{Sz} and
that in the case of rank two ($q=2$), the Jacobi polynomials of type BC were f\/irst studied~by
Koornwinder~\cite{K1,K2}.

}

\section[Random walks on the dual of a~compact Grassmannian and on $P_+$]{Random walks on the dual of a~compact Grassmannian\\ and on $\boldsymbol{P_+}$}

Recall that for integers $p\ge q$ the functions $\varphi_\lambda:=\varphi_\lambda^p$, $\lambda\in P_+$ form the spherical
functions of the compact Grassmannians $G/K=\mathcal G_{p,q}(\mathbb F)$.
As functions on~$G$, they are positive-def\/inite.
In other words, the Jacobi polynomials $(R_\lambda^p)_{\lambda\in P_+}$ are just the hypergroup characters of the
compact double coset hypergroups $G//K\cong A_0$.
We now recapitulate the associated dual hypergroup structures on~$P_+$.

\subsection[Dual hypergroup structures on $P_+$]{Dual hypergroup structures on $\boldsymbol{P_+}$}

As mentioned in the introduction, there is a~one-to-one correspondence between the set of (positive def\/inite) spherical
functions of $G/K$, which is parametrized by $P_+$, and the spherical unitary dual of~$G/K$, i.e.,~the set $\widehat G_K$
of all equivalence classes of irreducible unitary representa\-tions~$(\pi,H)$ of~$G$ whose restriction to~$K$ contains the
trivial representation with multiplicity one.
Here a~representation $(\pi,H)\in\widehat G_K $ and its spherical function $\varphi_\pi$ are related~by
\begin{gather*}
\varphi_\pi(x)=\langle u,\pi(x)u \rangle
\qquad
\text{for}
\quad
x\in G
\end{gather*}
with some~$K$-invariant vector $u\in H$ with $\|u\|=1$, which is determined uniquely up to a~complex constant of
absolute value 1.

Now consider $\lambda, \mu\in P_+$ with associated spherical functions $\varphi_\lambda$, $\varphi_\mu$ and the
corresponding representations $(\pi_\lambda,H_\lambda),(\pi_\mu,H_\mu)\in\widehat G_K $ with~$K$-invariant vectors
$u_\lambda$, $u_\mu$.
The tensor product $(\pi_\lambda\otimes\pi_\mu, H_\lambda\otimes H_\mu)$ is a~f\/inite-dimensional unitary representation
of~$G$ which decomposes into a~f\/inite orthogonal sum
\begin{gather*}
(\oplus_k\tau_k=\pi_\lambda\otimes\pi_\mu, \> \oplus_k \widetilde H_k=H_\lambda\otimes H_\mu)
\end{gather*}
of irreducible unitary representations $(\tau_k,\widetilde H_k)$ where some of them may appear several times.
Consider the orthogonal projections $p_k: H_\lambda\otimes H_\mu\to \widetilde H_k$.
Then the vectors $p_k(u_\lambda\otimes u_\mu)\in \widetilde H_k $ are~$K$-invariant, and for $p_k(u_\lambda\otimes
u_\mu)\ne0$, we obtain $(\tau_k,\widetilde H_k)\in \widehat G_K $, i.e., $(\tau_k,\widetilde H_k)$ is equal to some
$(\pi_\tau,H_\tau)$, $\tau\in P_+$.
Moreover, for $g\in G$,
\begin{gather*}
\varphi_\lambda(g)\varphi_\mu(g)=\langle u_\lambda\otimes u_\mu,\> (\pi_\lambda\otimes\pi_\mu)(g) u_\lambda\otimes u_\mu \rangle
\\
\phantom{\varphi_\lambda(g)\varphi_\mu(g)}
 =\sum\limits_k \langle p_k(u_\lambda\otimes u_\mu),\> \tau_k(g)\> p_k(
 u_\lambda\otimes u_\mu) \rangle
 =\sum\limits_k\|p_k(u_\lambda\otimes u_\mu)\|^2 \varphi_{\tau_k}(g)
\end{gather*}
with $\sum\limits_k\|p_k(u_\lambda\otimes u_\mu)\|^2=1$.
For $\lambda,\mu,\tau\in P_+$ we now def\/ine $c_{\lambda,\mu,\tau}\ge0$ as $\|p_k(u_\lambda\otimes u_\mu)\|^2$, whenever
$(\tau_k,\widetilde H_k)=(\pi_\tau,H_\tau)$ appears above with a~positive part, and $0$ otherwise.
As $\varphi_\lambda(g)\in\mathbb R$ for all~$\lambda$,~$g$ in our case, these nonnegative linearization coef\/f\/icients also
satisfy
\begin{gather*}
c_{\lambda,\mu,\tau}= \dim  H_\tau  \int_G \varphi_\lambda(g)\varphi_\mu(g)\varphi_\tau(g)\> dg.
\end{gather*}
For $\lambda,\mu\in P_+$ we def\/ine the probability measure
\begin{gather}
\label{def-Faltung}
\delta_\lambda*_{d,p}\delta_\mu:=\sum\limits_{\tau\in P_+}c_{\lambda,\mu,\tau}\delta_\tau\in M^1(P_+)
\end{gather}
with f\/inite support.
By its very construction, this convolution can be extended uniquely in a~weakly continuous, bilinear way to
a~probability preserving, commutative, and associative convolution on the Banach space $M_b(P_+)$ of all bounded, signed
measures on $P_+$.
Moreover, as all spherical functions are $\mathbb R$-valued in our specif\/ic examples, the contragredient representation
of any element in $\widehat G_K$ is the same representation, i.e., the canonical involution $.^*$ on $P_+\cong \widehat
G_K$ is the identity.
In summary, $(M_b(P_+), *_{d,p})$ is a~commutative Banach-$*$-algebra with the complex conjugation
$\mu^*(A):=\overline{\mu(A)}$ as involution.
Moreover, $(P_+,*_{d,p})$ becomes a~so-called Hermitian hypergroup in the sense of Dunkl, Jewett and Spector;
see~\cite{BH,Du,J}.

The Haar measure on this hypergroup, which is unique up to a~multiplicative constant, is the positive measure
$\omega=\sum\limits_{\lambda\in P_+} h(\lambda)\delta_\lambda$ with
\begin{gather*}
h(\lambda)= c_{\lambda,\lambda,0}^{-1}=\int_{G}\varphi_\lambda^2(g) \> dg =\frac{1}{\dim  H_\lambda},
\end{gather*}
where the f\/irst two equations follow from general hypergroup theory (see~\cite{J}) and the last one from the theory of
Gelfand pairs (see, e.g.,~\cite{F}).

The coef\/f\/icients $c_{\lambda,\mu,\tau}$ of the convolution $*_{d,p}$ on $P_+$ are related to the unique product
li\-nearization
\begin{gather*}
R_\lambda^p R_\mu^p=\sum\limits_{\tau\in P_+}c_{\lambda,\mu,\tau} R_\tau^p
\end{gather*}
of the Jacobi polynomials $R_\lambda^p$.
It is clear by our construction that for integers $p\ge q$, all $c_{\lambda,\mu,\tau}$ are nonnegative with
$\sum\limits_{\tau\in P_+}c_{\lambda,\mu,\tau}=1$.

Clearly, as $R_\lambda^p(0)=1$ for all~$\lambda$, the normalization also holds for all real $p\in[q,\infty[$.
We conjecture that actually all $c_{\lambda,\mu,\tau}$ are nonnegative for all $p\in[q,\infty[$ or at least for all
$p\in[2q-1,\infty[$.

Suppose that for f\/ixed $p\in[q,\infty[$ the linearization coef\/f\/icients $c_{\lambda,\mu,\tau}$ are all nonnegative.
Then equation~\eqref{def-Faltung} def\/ines a~commutative discrete hypergroup structure $(P_+,*_{d,p})$ with the
convolution
\begin{gather*}
\delta_\lambda*_{d,p}\delta_\mu:=\sum\limits_{\tau\in P_+}c_{\lambda,\mu,\tau}\delta_\tau\in M^1(P_+)
\end{gather*}
of point measures.
For instance, in the rank one case of Example~\ref{1-dim-bsp-1} the linearization coef\/f\/i\-cients~$c_{\lambda,\mu,\tau}$
are explicitly known and nonnegative for all $p\ge1$ as the product linearization coef\/f\/icients of the associated
one-dimensional Jacobi polynomials.

\subsection[Random walks on $P_+$]{Random walks on $\boldsymbol{P_+}$}\label{subsectionrandomwalk}

We next introduce certain random walks on $P_+$, i.e., time-homogeneous Markov chains on $P_+$ whose transition
probabilities are given in terms of the product linearizations coef\/f\/icients $c_{\lambda,\mu,\tau}$ for some f\/ixed
$p\in[q,\infty[$.
This concept even works, under a~suitable restriction, in the case where some of the $c_{\lambda,\mu,\tau}$ are negative.
To describe the restriction, we f\/ix $p\in[q,\infty[$ and def\/ine the set
\begin{gather*}
P_+^p:=\{\lambda\in P_+:\> c_{\lambda,\mu,\tau}\ge0
\;
\text{for all}
\;
\mu,\tau\in P_+\}
\end{gather*}
as well as
\begin{gather*}
M^1_p(P_+):=\big\{\nu\in M^1(P_+):\> \supp \nu\subset P_+^p\big\}.
\end{gather*}
We shall call the probability measures $\nu\in M^1_p(P_+)$ admissible.

Clearly, for integers $p\ge q$, as well as for $q=1$ we have $M^1_p(P_+)=M^1(P_+)$.
Unfortunately, it seems dif\/f\/icult to f\/ind further examples.
To illustrate the problem, consider the measure $\delta_{(1,0,\ldots,0)}\in M^1(P_+)$ for $q\ge 2$ and $p\in[q,\infty[$.
The explicit Pieri-type formula (6.4) of~\cite{D} then leads to a~concrete product linearization formula for
$R_\lambda^p  R_{(1,0,\ldots,0)}^p$.
It can be easily derived from~\cite{D} that for all $q\ge 2$, $p\in[q,\infty[$, and all $\lambda\ne \mu$, we have
$c_{\lambda,(1,0,\ldots,0), \mu}\ge0$ as desired.
However, we are so far not able to check from~\cite{D} that $c_{\lambda,(1,0,\ldots,0), \lambda}\ge0$ holds, which would
be necessary for $\delta_{(1,0,\ldots,0)}\in M^1_p(P_+)$.

As before, we f\/ix $d=1,2,4$ and $p\in [q,\infty[$.
We also f\/ix an admissible probability measure $\nu=\sum\limits_\mu p_\mu\delta_\mu\in M^1_p(P_+)$ and consider
a~time-homogeneous Markov chain $(S_n^{d,p})_{n\ge0}$ in discrete time on $P_+$ starting at time $0$ in the origin $0\in
P_+$ and with transition probability
\begin{gather*}
P\big(S_{n+1}^{d,p}=\tau|\> S_n^{d,p}
=\lambda\big)=(\nu*_{d,p}\delta_\lambda)(\{\tau\}),
\qquad
\lambda,\tau\in P_+,
\quad
n\in\mathbb N,
\end{gather*}
where
\begin{gather*}
\nu*_{d,p}\delta_\lambda:= \sum\limits_\tau \left(\sum\limits_\mu p_\mu c_{\lambda,\mu,\tau}\right) \delta_\tau\in M^1(P_+).
\end{gather*}
Such Markov-chains are called random walks on $(P_+,*_{d,p})$ associated with~$\nu$.
It is well-known and can be easily checked by induction that for all~$n$ the~$n$-fold convolution power
$\nu^{(n)}:=\nu*_{d,p}\cdots *_{d,p}\nu\in M^1(P_+)$ exists, and that $\nu^{(n)}$ is just the distribution
$P_{S_n^{d,p}}$ of $S_n^{d,p}$.

In view of the Central Limit Theorem~\ref{clt-spezial}, we give an interpretation of these convolution powers
$\delta_\lambda^{(n)}$ for integers $p\ge q$ and $\lambda\in P_+$ with $\lambda\ne 0$ in terms of representation theory.
We expect that this result is well-known, but we do not know an explicit proof in the literature.

\begin{Lemma}
\label{potenzen-tensoren}
Let $(\pi_\lambda,H_\lambda)$ be the non-trivial irreducible unitary representation of~$G$ associated with $\lambda\in
P_+$, $\lambda\ne0$ and with a~$K$-invariant vector $u_\lambda\in H_\lambda$ with $\|u_\lambda\|=1$.
For each $n\in \mathbb N$, decompose the~$n$-th tensor power $(\pi_\lambda^{\otimes, n},H_\lambda^{\otimes, n})$ into
its irreducible components
\begin{gather}
\label{tensor-zerl}
\big(\pi_\lambda^{\otimes, n},H_\lambda^{\otimes, n}\big)=\bigg(\bigoplus_{\tau_n}\pi_{\tau_n}, \bigoplus_{\tau_n} H_{\tau_n}\bigg)
\end{gather}
and consider the orthogonal projections $p_{\tau_n}:H_\lambda^{\otimes, n}\to H_{\tau_n}$.
Then for all $n\in\mathbb N$,
\begin{gather*}
\delta_\lambda^{(n)}= \sum\limits_{\tau_n}\big\|p_{\tau_n}\big(u_\lambda^{\otimes, n}\big)\big\|^2 \delta_{\tau_n}.
\end{gather*}
\end{Lemma}

\begin{proof}
We proceed by induction.
In fact, the case $n=1$ is trivial.
For $n\to n+1$, we start with~\eqref{tensor-zerl} and the associated orthogonal projections
$p_{\tau_n}:H_\lambda^{\otimes, n}\to H_{\tau_n}$.
We now decompose the products $\pi_{\tau_n}\otimes \pi_\lambda$ and obtain
\begin{gather*}
\big(\pi_\lambda^{\otimes, n+1},H_\lambda^{\otimes, n+1}\big)
=\bigg(\bigoplus_{\tau_n}(\pi_{\tau_n}\otimes \pi_\lambda),\bigoplus_{\tau_n} (H_{\tau_n}\otimes H_\lambda)\bigg)\\
\hphantom{\big(\pi_\lambda^{\otimes, n+1},H_\lambda^{\otimes, n+1}\big)}{}
=\bigg(\bigoplus_{\tau_n}\bigg(\bigoplus_{\mu_{k,n}}\pi_{\mu_{k,n}}\bigg),
\bigoplus_{\tau_n}\bigg(\bigoplus_{\mu_{k,n}}H_{\mu_{k,n}}\bigg)\bigg).
\end{gather*}
Notice that here the sum $\bigoplus_{\tau_n}\bigoplus_{\mu_{k,n}}$ corresponds to the sum $\bigoplus_{\tau_{n+1}}$ of
the lemma with $n+1$ instead of~$n$.
We now consider the orthogonal projections $p_{\mu_{k,n}}:H_\lambda^{\otimes, n+1}\to H_{\mu_{k,n}}$.
Then
\begin{gather*}
p_{\tau_{n}}\big(u_\lambda^{\otimes, n}\big)=c\big\| p_{\tau_{n}}\big(u_\lambda^{\otimes, n}\big)\big\|  u_{\tau_{n}},
\end{gather*}
where $|c|=1$, and thus
\begin{gather*}
\big\|p_{\mu_{k,n}}\big(u_\lambda^{\otimes, n+1}\big)\big\|^2  =\big\|p_{\mu_{k,n}}\big(p_{\tau_{n}}\big(u_\lambda^{\otimes, n}\big)\otimes u_\lambda\big)\big\|^2
 =\big\|p_{\tau_{n}}\big(u_\lambda^{\otimes, n}\big)\big\|^2 \big\|p_{\mu_{k,n}}\big(u_{\tau_{n}}\otimes u_\lambda\big)\big\|^2.
\end{gather*}
This fact, the assumption of our induction and the def\/inition of the convolution now readily imply the assertion of the
lemma for~$n+1$.
\end{proof}

We shall prove below that under a~natural moment condition on a~probability measure $\nu\in M^1(P_+)$, the~$C$-valued
random variables $\frac{1}{\sqrt n}S_n^{d,p}$ converge in distribution for $n\to\infty$.
In order to identify the limit distribution $\mu=\mu(d,p,\nu)\in M^1(C)$ in this central limit theorem, we need some
further preparations.

\subsection[Bessel convolutions on $C$ and Laguerre ensembles]{Bessel convolutions on $\boldsymbol{C}$ and Laguerre ensembles}

As described in Section~\ref{Section2}, the Bessel functions $\widetilde\varphi_\lambda^p$ with $\lambda\in C $ make up the set of
bounded spherical functions of the Euclidean symmetric space $G_0/K$ with $K=U(p,\mathbb F)\times U(q,\mathbb F)$ and
$G_0= K\ltimes M_{p,q}(\mathbb F)$.
Thus by the notation of~\cite{BH} and~\cite{J}, the chamber $C\cong G_0//K$ with the associated double coset convolution
$\bullet_{d,p}$ is a~commutative double coset hypergroup with the functions $\widetilde\varphi_\lambda^p $ as (bounded)
hypergroup characters.
We now introduce the probability measures
\begin{gather*}
d\rho_{d,p}(x)
:=c_{d,p}^{-1} \prod\limits_{j=1}^q x_j^{d(p-q+1)-1}  \prod\limits_{1\le i<j\le q}\big(x_i^2-x_j^2\big)^d  e^{-(x_1^2+\cdots+x_q^2)/2}\> dx
\end{gather*}
on the Weyl chamber~$C$, with the normalization constant
\begin{gather*}
c_{d,p}=\int_{C} \prod\limits_{j=1}^q x_j^{d(p-q+1)-1}   \prod\limits_{1\le i<j\le q}\big(x_i^2-x_j^2\big)^d  e^{-(x_1^2+\cdots+x_q^2)/2}\> dx.
\end{gather*}
The measure $\rho_{d,p}\in M^1(C)$ is well-known in the random matrix theory of so-called Laguerre- or~$\beta$-ensembles
as it is the distribution of the singular values of a~$M_{p,q}(\mathbb F)$-valued random variable for which the real and
imaginary parts of all entries are i.i.d.\
and standard normally distributed.
This fact is well-known; it can also be derived from the Haar measure of the double coset hypergroups
$(C,\bullet_{d,p})$ in~\cite[Section~4.1]{R1}.
Moreover, having this group theoretic interpretation in mind, one easily obtains the following well-known relation from
the Fourier transform of a~multivariate standard normal distribution:
\begin{gather}
\label{Dunkl-Gaussian-eq}
\int_{C} \widetilde\varphi_\lambda^p(x)\> d\rho_{d,p}(x)= e^{-(\lambda_1^2+\cdots+\lambda_1^2)/2}
\qquad
\text{for}
\quad
\lambda\in C,
\end{gather}
see~\cite[Proposition~XV.2.1]{FK} or~\cite{V0}.
This identif\/ication of the spherical Fourier-transform of~$\rho_{d,p}$ will be essential in the following for the
central limit theorem.

\subsection[A central limit theorem for random walks on $P_+$]{A central limit theorem for random walks on $\boldsymbol{P_+}$}

The probability measure $\rho_{d,p}$ appears in the CLT below as limit up to some scaling parameter
$\sigma^2=\sigma^2(\nu,p,d)$, which admits an interpretation as a~variance parameter.
For the description of $\sigma^2$, we need the so-called moment functions on $(P_+,*_{d,p})$ up to order two.
For the general theory of such moment functions and their applications to limit theorems for random walks on hypergroups
we refer to \cite[Chapter~7]{BH} and \cite{Z}, and the references there.
The moment functions are characterized by additive functional equations similar to the multiplicative ones for
hypergroup characters.
They are usually def\/ined in terms of (partial) derivatives with respect to the spectral variables at the identity
character.
In our examples, the identity corresponds to $x=0\in C$.
This motivates the following def\/inition.

\begin{Definition}
Let $p\in[q,\infty[$ be f\/ixed.
For $k,l=1,\ldots,q$ we def\/ine the moment functions $m_k, m_{k,l}: P_+\to\mathbb R$ of the f\/irst and second order~by
\begin{gather*}
m_k(\lambda):= i  \frac{\partial}{\partial x_k} R_\lambda(x)\bigg|_{x=0}
\qquad
\text{and}
\qquad
m_{k,l}(\lambda):= -\frac{\partial^2}{\partial x_k\partial x_l} R_\lambda(x)\bigg|_{x=0}.
\end{gather*}
\end{Definition}

\noindent Let us collect some properties of these moment functions.

\begin{Lemma}\label{basic-prop-moment}\quad
\begin{enumerate}\itemsep=0pt
\item[$(1)$]
For all $k,l=1,\ldots,q$ with $k\ne l$ and all $\lambda\in P_+$, $m_k(\lambda)= m_{k,l}(\lambda)=0$.
\item[$(2)$]
The functions $m_{k,k}$ are independent of $k=1,\ldots,q$, and the function $m:=m_{1,1}: P_+\to\mathbb
R$ is a~second order polynomial of the form
\begin{gather*}
m(\lambda)
=\frac{1}{4}\sum\limits_{r,s=1}^q a_{r,s}
 (\lambda_r-\lambda_{r+1})(\lambda_s-\lambda_{s+1}) +\frac{1}{2} \sum\limits_{r=1}^q b_r (\lambda_r-\lambda_{r+1})
\end{gather*}
with suitable coefficients $a_{r,s}$, $b_r$.
In particular, $m(0)=0$.
\item[$(3)$]
For all $\lambda,\tau\in P_+$, $\int_{P_+} m\> d(\delta_\lambda*_{d,p}\delta_\tau)= m(\lambda)+m(\tau)$.
\end{enumerate}
\end{Lemma}

\begin{proof}
The Jacobi polynomials $R_\lambda(x)$ are invariant under the Weyl group of type B acting in the variable~$x$.
In particular, $R_\lambda(x_1,\ldots,x_q)$ is even in each $x_i$, and this gives part (1).
Moreover, as $R_\lambda(x_1,\ldots,x_q)$ is invariant under permutations of the $x_i$, the function $m_{k,k}$ is
independent of~$k$.
We now study $m:=m_{1,1}$ more closely.
We start with the case $p> 2q-1$.
In this case we obtain from the integral representation~\eqref{int-rep} that
\begin{gather}
\label{m-int-rep}
m(\lambda)= -\int_{B_q\times U(q,\mathbb F)} \frac{\partial^2}{\partial x_1^2}\left(\Delta_{\lambda/2}(g_{ix}(u,w))\right)\bigg|_{x=0}\> dm_p(w)du
\end{gather}
with the power function
\begin{gather*}
\Delta_{\lambda/2}(g_{ix}(u,w)) =\prod\limits_{r=1}^q \Delta_r(g_{ix}(u,w))^{(\lambda_r-\lambda_{r+1})/2}
\qquad
\text{with}
\quad
\lambda_{q+1}=0.
\end{gather*}
A~short calculation, using that $\Delta_{\lambda/2}(g_0(u,w))=1$, gives
\begin{gather}
\frac{\partial^2}{\partial x_1^2} \Delta_{\lambda/2}(g_{ix}(u,w))\bigg|_{x=0}
=\frac{1}{4} \left(\sum\limits_{r=1}^q\frac{\partial}{\partial x_1}\left(\Delta_r(g_{ix}(u,w))\right)\bigg|_{x=0}
 (\lambda_r-\lambda_{r+1})\right)^2
\nonumber
\\
\hphantom{\frac{\partial^2}{\partial x_1^2} \Delta_{\lambda/2}(g_{ix}(u,w))\bigg|_{x=0}=}{}
+\frac{1}{2}\sum\limits_{r=1}^q \frac{\partial^2}{\partial x_1^2} \left(\ln \Delta_r(g_{ix}(u,w))\right)\bigg|_{x=0}
  (\lambda_r-\lambda_{r+1}).
\label{ableitung-power}
\end{gather}
Formulas~\eqref{m-int-rep} and~\eqref{ableitung-power} now imply that~$m$ is a~second order polynomial as claimed, with
linear terms
\begin{gather*}
b_r= \int_{B_q\times U(q,\mathbb F)} \frac{\partial^2}{\partial x_1^2} \left(\ln\Delta_r(g_{ix}(u,w))\right)\bigg|_{x=0}\> dm_p(w)du.
\end{gather*}
The coef\/f\/icients $a_{r,s}$ are obtained from the Taylor expansion~\eqref{entwickel} of $\Delta_r(g_{ix}(u,w))$. They are
given~by
\begin{gather}
\label{a_rs}
a_{r,s}:= \int_{B_q\times U(q,\mathbb F)} \operatorname{tr}_r(u^{*}(w^*P_1+ P_1  w)u) \operatorname{tr}_s(u^{*}(w^*P_1+ P_1  w)u)\>dm_p(w)du
\end{gather}
with the diagonal matrix $P_1:=\diag(1,0,\ldots,0)\in M_q(\mathbb F)$.
This proves that~$m$ is a~second order polynomial for $p>2q-1$.
The case $p\ge q$ follows by analytic continuation.

Finally, for the proof of part (3) we observe that for $\lambda,\tau\in P_+$,
\begin{gather*}
\int_{P_+} m\> d(\delta_\lambda*_{d,p}\delta_\tau)
=-\int_{P_+} \frac{\partial^2}{\partial x_1^2}R_\kappa(x)\bigg|_{x=0} \> d(\delta_\lambda*_{d,p}\delta_\tau)(\kappa)
\\
\phantom{\int_{P_+} m\> d(\delta_\lambda*_{d,p}\delta_\tau)}
=-\frac{\partial^2}{\partial x_1^2}\left(\int_{P_+} R_\kappa(x)\>d(\delta_\lambda*_{d,p}\delta_\tau)\right)\bigg|_{x=0}
\\
\phantom{\int_{P_+} m\> d(\delta_\lambda*_{d,p}\delta_\tau)}
=-\frac{\partial^2}{\partial x_1^2}\left(R_\lambda(x)R_\tau(x)\right)\bigg|_{x=0} \>=\> m(\lambda)+m(\tau).
\end{gather*}
Notice that the last equality follows from part (1) and $R_\lambda(0)=1$.
\end{proof}

\begin{Example}
[the rank one case] For $q=1$, the moment function~$m$ is given~by
\begin{gather*}
m(\lambda)= \frac{\lambda(\lambda+dp+d-2)}{dp},
\qquad
\lambda\in 2\mathbb Z^+.
\end{gather*}
In fact, this can be easily derived from the def\/inition of~$m$ and the explicit formulas for the classical Jacobi
polynomials in~\eqref{class-Jacobi1} and~\eqref{class-Jacobi2}.
Moreover, it can be also derived from the proof of part (2) of the preceding lemma and a~direct elementary computation
of $a_{1,1}$ and $b_1$ for $q=1$ there.
\end{Example}

We shall need the following variant of Lemma~\ref{basic-prop-moment}(2) on the growth of~$m$ for $p\in[q,\infty[$:

\begin{Lemma}
\label{quad-wachstum-oben}
For all $x\in\mathbb R^q$, $\lambda\in P_+$ and $k,l=1,\ldots,q$,
\begin{gather*}
\left|\frac{\partial}{\partial x_k} R_\lambda^p(x)\right|\leq \lambda_1
\qquad
\text{and}
\qquad
\left|\frac{\partial^2}{\partial x_k\partial x_l} R_\lambda^p(x)\right|\leq \lambda_1^2.
\end{gather*}
In particular, $m(\lambda)\leq \lambda_1^2$ for $\lambda\in P_+$.
\end{Lemma}

\begin{proof}
Let again~$W$ denote the Weyl group of type $BC_q$.
We introduce the normalized~$W$-invariant orbit sums
\begin{gather*}
\widetilde M_\lambda (x):= \frac{1}{|W\lambda|} \sum\limits_{\mu\in W.\lambda}e^{i\langle \lambda, x\rangle},
\qquad
\lambda \in P_+.
\end{gather*}
Then the Jacobi polynomials $R_\lambda^p$ can be written as linear combinations of such orbit sums.
It follows from the considerations in~\cite[Section~11]{M} that for non-negative multiplicity values, the expansion
coef\/f\/icients are all non-negative.
That is,
\begin{gather*}
R_\lambda^p=\sum\limits_{\mu\in P_+: \mu \leq \lambda} c_{\lambda\mu} \widetilde M_\mu
\end{gather*}
with
\begin{gather*}
c_{\lambda\mu}=c_{\lambda\mu}^p\geq 0, \qquad \sum\limits_{\mu \leq \lambda} c_{\lambda\mu}=1.
\end{gather*}
Here $\leq$ denotes the dominance order on $P_+$ given by $\mu \leq \lambda
\Longleftrightarrow
\sum\limits_{i=1}^r \mu_i \leq \sum\limits_{i=1}^r \lambda_i $ for $ r= 1, \ldots, q$. We have
\begin{gather*}
\partial_{x_k} \widetilde M_\lambda(x)=\frac{1}{|W\lambda|} \sum\limits_{\mu\in W.\lambda} i\mu_k e^{i\langle\mu,x\rangle},
\qquad
\partial_{x_kx_l} \widetilde M_\lambda(x)=\frac{-1}{|W\lambda|} \sum\limits_{\mu\in W.\lambda} \mu_k\mu_l e^{i\langle\mu,x\rangle}.
\end{gather*}
Notice that $|\mu_k| \leq \lambda_k \leq \lambda_1$ for each $\mu\in W\lambda$.
We thus obtain, independently of $x\in A_0$,
\begin{gather*}
\big\vert \partial_{x_k} \widetilde M_\lambda(x)\big\vert \leq \lambda_1,
\qquad
\big\vert \partial_{x_k}\partial_{x_l} \widetilde M_\lambda(x)\big\vert \leq \lambda_1^2.
\end{gather*}
Further, if $\mu\in P_+$ with $\mu \leq \lambda$, then $\mu_1\leq \lambda_1$ and therefore
\begin{gather*}
\big\vert \partial_{x_k} R_\lambda^p(x)\big\vert \leq \sum\limits_{\mu \leq \lambda} c_{\lambda\mu} \big\vert \partial_{x_k}
\widetilde M_\mu(x)\big\vert \leq \sum\limits_{\mu \leq \lambda}c_{\lambda\mu} \mu_1 \leq \lambda_1.
\end{gather*}
In the same way,
\begin{gather*}
\big\vert \partial_{x_k}\partial_{x_l} R_\lambda^p(x)\big\vert \leq \sum\limits_{\mu \leq \lambda} c_{\lambda\mu}\big\vert
\partial_{x_k}\partial_{x_l} \widetilde M_\mu(x)\big\vert \leq \sum\limits_{\mu \leq \lambda}c_{\lambda\mu}\mu_1^2\leq\lambda_1^2.\tag*{\qed}
\end{gather*}
\renewcommand{\qed}{}
\end{proof}

We also need some further properties of the moment function~$m$.
We here have to restrict our attention to the case $p\in \{q,q+1,\ldots,2q-1\}\cup]2q-1,\infty[$.
We shall assume this restriction from now on.
We expect that the results below are also valid for all $p\in[q,\infty[$.

\begin{Lemma}\label{quad-wachstum}\qquad
\begin{enumerate}\itemsep=0pt
\item[$(1)$] The matrix $A=(a_{r,s})_{r,s=1,\ldots,q}\in M_q(\mathbb R)$ is positive definite.
\item[$(2)$] For all $\lambda\in P_+\setminus\{0\}$, $m(\lambda)>0$.
\item[$(3)$] There exists a~constant $C_1>0$ such that for all $\lambda\in P_+$, $C_1\lambda_1^2\le m(\lambda)$.
\end{enumerate}
\end{Lemma}

\begin{proof}
For the proof of (1), we f\/irst consider the case $p\in{}]2q-1,\infty[$ and conclude from the def\/inition of the $a_{r,s}$
in the proof in Lemma~\ref{basic-prop-moment} that~$A$ is symmetric, and that for all $\tau\in \mathbb R^q$,
\begin{gather*}
\tau^T A\tau= \int_{B_q\times U(q,\mathbb F)} \left(\sum\limits_{r=1}^q \tau_r  \operatorname{tr}_r(\cdot u^{*}(w^*P_1+ P_1w)u)\right)^2\> dm_p(w)du,
\end{gather*}
where the functions
\begin{gather*}
(w,u)\mapsto \operatorname{tr}_r(\cdot u^{*}(w^*P_1+ P_1 w)u),
\qquad
B_q\times U(q,\mathbb F)\to \mathbb R,
\end{gather*}
are linearly independent for $r=1,\ldots q$.
This shows that $\tau^T A\tau>0$ for all $\tau\in \mathbb R^q$ with $\tau\ne 0$ as claimed.
The case of integers $p\ge q$ can be handled in a~similar way by using a~modif\/ied version of integral
representation~\eqref{a_rs} for $a_{r,s}$ which is based on Remark~\ref{group} instead of
Proposition~\ref{int-rep-prop}.

For the proof of part (2) we proceed as in the proof of Lemma~\ref{quad-wachstum-oben} and write
\begin{gather*}
R_\lambda^p=\sum\limits_{\mu\in P_+: \mu \leq \lambda} c_{\lambda\mu} \widetilde M_\mu
\qquad
\text{with}
\quad
c_{\lambda\mu} \ge 0, \quad c_{\lambda\lambda}>0,
\quad
\sum\limits_{\mu \leq \lambda} c_{\lambda\mu}=1.
\end{gather*}
Thus for $\lambda\in P_+\setminus\{0\}$,
\begin{gather*}
m(\lambda)=-\partial_{x_1}^2 R_\lambda^p(x)\big|_{x=0}= \sum\limits_{\mu\in P_+: \mu \leq \lambda}
\frac{c_{\lambda\mu}}{|W\mu|}\sum\limits_{\tau\in W\mu} \tau_1^2>0.
\end{gather*}

For the proof of part (3), we use (1) and write $m(\lambda)$ as
\begin{gather*}
m(\lambda)= \lambda^T \widetilde A\lambda - b^T\lambda
\end{gather*}
with some positive def\/inite matrix $\widetilde A$ and some $b\in\mathbb R^q$.
We thus f\/ind constants $c,d>0$ such that $m(\lambda)-c\lambda_1^2>0$ holds for all $\lambda\in P_+$ with $\lambda_1\ge
d$.
As there are only f\/initely many $\lambda\in P_+$ with $\lambda_1< d$, we conclude from part (2) that there exists some
$C_1>0$ with $m(\lambda)-c\lambda_1^2>0$ for all $\lambda\in P_+$ with $\lambda\ne 0$.
\end{proof}

\begin{Remark}
The nonnegativity of $m(\lambda)$ in Lemma~\ref{quad-wachstum}(2) can be easily established directly.
In fact, assume that $m(\lambda)<0$ for some $\lambda\in P_+$.
Then the Taylor formula
\begin{gather*}
R_\lambda(x)=1-\frac{m(\lambda)}{2}\big(x_1^2+\cdots+x_q^2\big)+O\big(\|x\|^3\big)
\end{gather*}
implies that $R_\lambda(x)>1$ for some~$x$ close to~$0$, and thus by the Weyl group invariance of $R_\lambda$, for some
$x\in A_0$.
But this contradicts the fact that $\|R_\lambda\|_\infty\le1$ on $A_0$, which is a~consequence of
Proposition~\ref{int-rep-prop}.
However, we have no dif\/ferent proof for the strict positivity of $m(\lambda)$ for $\lambda\ne 0$ than the one given in
Lemma~\ref{quad-wachstum}(2).
\end{Remark}

We next use the moment function~$m$ in order to def\/ine a~modif\/ied variance of measures $\nu\in M^1(P_+)$ depending on
the underlying convolution $*_{d,p}$.
This modif\/ied variance will appear in the CLT below.

\begin{Definition}
%\label{def-variance}
Let $\nu\in M^1(P_+)$ be a~probability measure with f\/inite second moments, meaning that $\sum\limits_{\lambda\in P_+}
\lambda_1^2 \nu(\{\lambda\})<\infty$. Then the modif\/ied second moment $\sigma^2:=\sigma^2(\nu)$ of~$\nu$ is def\/ined as
\begin{gather*}
\sigma^2:=\sum\limits_{\lambda\in P_+} m(\lambda) \nu(\{\lambda\}).
\end{gather*}
Notice that by Lemmas~\ref{basic-prop-moment} and~\ref{quad-wachstum}, $\sigma^2$ is f\/inite and non-negative where
$\sigma^2=0$ holds precisely for $\nu=\delta_0$.
\end{Definition}

\begin{Lemma}
\label{diffbar}
Let $\nu\in M^1(P_+)$ be a~probability measure with finite second moments and with the modified variance $\sigma^2\ge0$.
Then the spherical Fourier transform
\begin{gather*}
{\cal F}\nu: \ \mathbb R^q\to\mathbb R
\qquad
\text{with}
\quad
{\cal F}\nu(x):=\sum\limits_{\lambda\in P_+} R_\lambda^p(x) \nu(\{\lambda\})
\end{gather*}
is twice continuously differentiable on $\mathbb R^q$ with
\begin{gather*}
{\cal F}\nu(0)=1,
\qquad
\nabla{\cal F}\nu(0)=0,
\qquad
\text{and Hesse matrix}
\quad
D^2({\cal F}\nu)(0)=-\sigma^2  I_q.
\end{gather*}
\end{Lemma}

\begin{proof}
As~$\nu$ has f\/inite second moments, we conclude form Lemma~\ref{quad-wachstum-oben} that for all $k,l=1,\ldots,q$, the
series
\begin{gather*}
\sum\limits_{\lambda\in P_+} \nu( \{\lambda\})
 \frac{\partial}{\partial x_k} R_\lambda^p(x),
\qquad
\sum\limits_{\lambda\in P_+} \nu(\{\lambda\})   \frac{\partial^2}{\partial x_k\partial x_l} R_\lambda^p(x)
\end{gather*}
converge uniformly with respect to $x\in\mathbb R^q$.
This implies that ${\cal F}\nu(x)$ is twice continuously dif\/ferentiable on $\mathbb R^q$, and the partial derivatives
commute with the summation.
The derivatives at $x=0$ are now obtained by Lemma~\ref{basic-prop-moment}.
\end{proof}

We now turn to random walks $(S_n^{d,p})_{n\ge0}$ on $(P_+,*_{d,p})$ associated with some admissible $\nu\in
M^1_p(P_+)$.
It is well-known (see, e.g.,~\cite[Section 7.3]{BH}) that the additive functional equation for~$m$ in
Lemma~\ref{basic-prop-moment}(3) leads to relations between the modif\/ied variance of~$\nu$ and random walks associated
with~$\nu$.

\begin{Lemma}%\label{addi-moments}
\qquad
\begin{enumerate}\itemsep=-1pt
\item[$(1)$]
For all $\nu_1,\nu_2\in M^1_p(P_+)$ with finite second moments, the measure $\nu_1*_{d,p}\nu_2\in
M^1(P_+)$ has also finite second moments, and
$\sigma^2(\nu_1*_{d,p}\nu_2)= \sigma^2(\nu_1)+\sigma^2(\nu_2)$.

\item[$(2)$]
Let $(S_n^{d,p})_{n\ge0}$ be a~random walk on $(P_+,*_{d,p})$ associated with the measure $\nu \in
M^1_p(P_+)$ with finite second moments.
Then, for all integers $n\ge0$, the expectation of $m(S_n^{d,p})$ satisfies $E(m(S_n^{d,p}))=n\sigma^2(\nu)$, and the
process $(m(S_n^{d,p})-n\sigma^2(\nu))_{n\ge0}$
is a~martingale with respect to the canonical filtration of $(S_n^{d,p})_{n\ge0}$.
\end{enumerate}
\end{Lemma}

\begin{proof}
Part (1) follows easily from Lemma~\ref{basic-prop-moment}(3); c.f.~\cite[Section~7.3.7]{BH}.
Moreover, the f\/irst assertion of~(2) follows from~(1) by induction.
For the proof of the second statement in~(2) we refer to~\cite[Proposition 7.3.19]{BH}.
\end{proof}

\begin{Lemma}
\label{sn-wachstum}
For $a\in{}]0,\infty[$, define the finite set $K_a:=\{\lambda\in P_+:\> \lambda_1\le a\}$.
Let $(S_n^{d,p})_{n\ge0}$ be a~random walk on $(P_+,*_{d,p})$ as described above associated with some admissible $\nu\in
M^1_p(P_+)$ with finite second moments.
Then, for each $\epsilon>0$ there exists some $a\ge1$ such that for all $n\in\mathbb N$, $P(S_n^{d,p}\not\in K_{\sqrt
n  a})\le\epsilon$.
\end{Lemma}

\begin{proof}
By Lemma~\ref{quad-wachstum}, we f\/ind $c>0$ with $m(\lambda)\ge c\lambda_1^2$ for $\lambda\in P_+$.
Therefore, $m(\lambda)\ge c\lambda_1^2\ge ca^2n$ for all $\lambda\in P_+\setminus K_{\sqrt n  a}$ and all~$n$.
Hence,
\begin{gather*}
P\big(S_n^{d,p}\not\in K_{\sqrt n  a}\big)  \le  P\big(m\big(S_n^{d,p}\big)\ge ca^2n\big)   \le   \frac{E(m(S_n^{d,p}))}{ca^2n}
 = \frac{\sigma^2n}{ca^2n}=\frac{\sigma^2}{ca^2}
\end{gather*}
with the f\/inite modif\/ied variance~$\sigma^2$.
This implies the claim.
\end{proof}

We are now ready to prove the main result of this section:

\begin{Theorem}\label{clt}\looseness=-1
Let $p\in \{q,q+1,\ldots,2q-1\}\cup]2q-1,\infty[$.
and $\nu\in M^1_p(P_+)$ be an admissible probability measure with $\nu\ne \delta_0$ and with finite second moments.
Let $\sigma^2\in{}]0,\infty[$ be the modified modified variance of~$\nu$, and $(S_n^{d,p})_{n\ge0}$ be a~random walk on
$(P_+,*_{d,p})$ associated with~$\nu$.
Then $ S_n^{d,p}/\sqrt{n\sigma^2}$ converges in distribution to the distribution $\rho_{d,p}\in M^1(C)$ of a~Laguerre
ensemble in~$C$.
\end{Theorem}

\begin{proof}
Fix $x\in C$.
Let $n\in\mathbb N$ be large enough such that $x/\sqrt n\in A_0$.
By Section~\ref{subsectionrandomwalk}, $S_n^{d,p}$ has the distribution $\nu^{(n)}$.
We thus obtain from the multiplicativity of the spherical Fourier transform of probability measures on $(P_+,*_{d,p})$,
Lemma~\ref{diffbar}, the qualitative Taylor formula, and from the properties of the moment functions in
Lemma~\ref{basic-prop-moment} that for $n\to \infty$,
\begin{gather}
E\big(R_{S_n^{d,p}}(x/\sqrt n)\big) =\big({\cal F}\nu(x/\sqrt n)\big)^n=\left(1-\frac{\sigma^2}{2n}\big(x_1^2+\cdots +x_q^2\big)+o(1/n)\right)^n
\nonumber\\
\hphantom{E\big(R_{S_n^{d,p}}(x/\sqrt n)\big)}{}
\longrightarrow e^{-(x_1^2+\cdots +x_q^2)\sigma^2/2}.\label{exp-limit}
\end{gather}

Now f\/ix $\epsilon>0$.
By Lemma~\ref{sn-wachstum}, there is so $a>0$ such that for all $n\in\mathbb N$, $P(S_n^{d,p}\not\in K_{\sqrt n
a})\le\epsilon$.
We now conclude from the Mehler--Heine formula~\eqref{vgl-Bessel} that for all $\lambda\in K_{\sqrt n a}$, that is,
$\lambda\in P_+$ with $\lambda_1\le a\sqrt n $,
\begin{gather*}
\big|R_{\lambda}^p(x/\sqrt n)-\widetilde\varphi_\lambda^{p}(x/\sqrt n)\big|\leq C_1  x_1^2\lambda_1  e^{C_2x_1^2\lambda_1/\sqrt n}\le\epsilon
\end{gather*}
whenever~$n$ is suf\/f\/iciently large.
As $|R_{\lambda}^p(x/\sqrt n)|\le1$ and $|\widetilde\varphi_\lambda^{p}(x/\sqrt n)|\le1$, we thus have
\begin{gather*}
\big|E \big(R_{S_n^{d,p}}(x/\sqrt n)\big) - E\big(\widetilde\varphi_{S_n^{d,p}}(x/\sqrt n)\big)\big|
\\
\qquad
 \le E\left(\big|R_{S_n^{d,p}}(x/\sqrt n)-\widetilde\varphi_{S_n^{d,p}}(x/\sqrt n)\big| {\bf 1}_{\{S_n^{d,p}\in
K_{\sqrt n a}\}} \right)+2  P\big(S_n^{d,p}\not\in K_{\sqrt n  a}\big)
 \le 3\epsilon
\end{gather*}
for~$n$ suf\/f\/iciently large.
Together with~\eqref{exp-limit} and the identity $\widetilde\varphi_{cy}^p(x)=\widetilde\varphi_{y}^p(cx)$ for $c>0$ and
$x,y\in C$, this implies that for all $x\in C$,
\begin{gather*}
\lim_{n\to\infty} E\big(\widetilde\varphi_{S_n^{d,p}/\sqrt{\sigma^2n}}^p(x)\big)
=\lim_{n\to\infty}E\big(\widetilde\varphi_{S_n^{d,p}}^p\big(x/\sqrt{\sigma^2n}\big)\big)
\\
\hphantom{\lim_{n\to\infty} E\big(\widetilde\varphi_{S_n^{d,p}/\sqrt{\sigma^2n}}^p(x)\big)}{}
 =\lim_{n\to\infty} E\big(R_{S_n^{d,p}}\big(x/\sqrt{\sigma^2n}\big)\big)=e^{-(x_1^2+\cdots+x_q^2)/2}.
\end{gather*}
From this limit, equation~\eqref{Dunkl-Gaussian-eq} and Levy's continuity theorem for the spherical Fourier transform on
the double coset hypergroup $(C, \bullet_{d,p})$ (see, e.g.,
\cite[Section~4.2]{BH}), we now infer that $S_n^{d,p}/\sqrt{\sigma^2n}$ converges in distribution to $\rho_{d,p}$ as
claimed.
\end{proof}

Theorem~\ref{clt-spezial} in the introduction is an immediate consequence from Theorem~\ref{clt} and
Lem\-ma~\ref{potenzen-tensoren}.

We also remark that the methods of the preceding proof lead with some additional technical ef\/fort to rates of
convergence in the CLT; see~\cite{Ga,V1} for the rank one case.

We f\/inish this paper with a~strong law of large numbers; it follows easily from the preceding properties of the moment
function~$m$, in combination with strong laws of large numbers for random walks on commutative hypergroups in~\cite[Section~7.3]{BH} and~\cite{Z}.

\begin{Theorem}
%\label{LLN}
Let $\nu\in M^1_p(P_+)$ be admissible with with finite second moments, and let $(S_n^{d,p})_{n\ge0}$ be an associated
random walk on $(P_+,*_{d,p})$.
Then for all $\epsilon>1/2$, $S_n/n^\epsilon\to0$ almost surely.
\end{Theorem}

\begin{proof}
Consider f\/irst the hypergroup case with an integer $p\ge q$.
By Lemmas~\ref{basic-prop-moment} and~\ref{quad-wachstum}, all conditions of Theorem 7.3.26 in~\cite{BH} are satisf\/ied
for the time-homogeneous random walk $(S_n^{d,p})_{n\ge0}$, the sequence $(r_n:= n^{2\epsilon})_{n\ge1}$, and the moment
function~$m$ instead of $m_2$ in~\cite{BH}.
This theorem now yields that $m(S_n^{d,p})/ n^{2\epsilon}$ tends to $0$ almost surely, and Lemma~\ref{quad-wachstum}
proves the claim.
An inspection of the details in the proof of Theorem 7.3.26 in~\cite{BH} shows that this theorem is also available for
all~$p$ and admissible $\nu\in M^1_p(P_+)$ which proves the theorem in general.
\end{proof}

\vspace{-4mm}

\pdfbookmark[1]{References}{ref}
\LastPageEnding

\end{document}